%% file: paper_hpc_julia.tex
\documentclass[Afour,sageh,times]{sagej}

\usepackage{moreverb,url}
\usepackage[colorlinks,bookmarksopen,bookmarksnumbered,citecolor=red,urlcolor=red]{hyperref}
\usepackage{siunitx}
\usepackage{xcolor}
\usepackage{tabularx}
\usepackage{cuted}

\newcommand\BibTeX{{\rmfamily B\kern-.05em \textsc{i\kern-.025em b}\kern-.08em
T\kern-.1667em\lower.7ex\hbox{E}\kern-.125emX}}

\def\volumeyear{2026}

\newcommand{\blue}[1]{{\color{blue} #1}}
\newcommand{\red}[1]{{\color{red} #1}}

\runninghead{Candelaresi et al.}

\title{Massively parallel numerical simulations with Julia}

\author{%
Simon Candelaresi\affilnum{1},
Benedict Geihe\affilnum{2},
Marco Artiano\affilnum{3},
Lars Christmann\affilnum{2, 4},
Valentin Churavy\affilnum{1, 3},
Andrés Rueda-Ramírez\affilnum{5},
Hendrik~Ranocha\affilnum{3},
Gregor~J.~Gassner\affilnum{3},
Michael~Schlottke-Lakemper\affilnum{1, 4}
}

\affiliation{
\affilnum{1}High-Performance Scientific Computing, Centre for Advanced Analytics and Predictive Sciences, University of Augsburg, Germany\\
\affilnum{2}Department of Mathematics and Computer Science, Division of Mathematics, University of Cologne, Germany\\
\affilnum{3}Institute of Mathematics, Johannes Gutenberg University Mainz, Germany\\
\affilnum{4}Applied and Computational Mathematics, RWTH Aachen University, Germany\\
\affilnum{5}School of Aerospace Engineering, Universidad Politécnica de Madrid, Spain
}

\corrauth{Simon Candelaresi, High-Performance Scientific Computing, Centre for Advanced Analytics and Predictive Sciences, University of Augsburg, Germany.}
\email{simon.candelaresi@proton.me}

\begin{document}

\begin{abstract}
  \input{abstract.tex}
\end{abstract}

\keywords{high-performance computing, parallel scalability, Julia, MPI}

\maketitle

\input{mainmatter.tex}

\section*{Acknowledgments}

\begin{acks}
\input{acknowledgments.tex}
\end{acks}

\bibliographystyle{SageH}
\bibliography{references}

\end{document}

%% file: abstract.tex
The Julia programming language aims to provide a modern approach to develop 
high-performance computing (HPC) applications. It tries to achieve this by combining a
high-level, dynamic interface with just-in-time compilation to native machine
code, thereby facilitating high developer productivity and native code
performance at the same time. While this approach has already been shown to work well
for serial applications, it is not clear if it readily translates to traditional,
massively parallel HPC work loads. In this paper, we fill this gap by
analyzing the parallel performance of the numerical computational fluid dynamics simulation code Trixi.jl,
written in Julia, and compare it to the Fortran code FLUXO.
We show some of the
challenges of using Julia at scale and discuss possible solutions, specifically
with respect to code loading and compilation at startup.
Finally, we demonstrate the parallel scaling of our Julia code on up to \num{61440} CPU cores.

%% file: mainmatter.tex
%%%%%%%%%%%%%%%%%%%%%%%%%%%%%%%%%%%%%%%%%%%%%%%%%%%%%%%%%%%%%%%%%%%%%%%%%%%%%%%%%%%%%%%%%%%%
%%%%%%%%%%%%%%%%%%%%%%%%%%%%%%%%%%%%%%%%%%%%%%%%%%%%%%%%%%%%%%%%%%%%%%%%%%%%%%%%%%%%%%%%%%%%
\section{Introduction}\label{sec:introduction}

Julia's\footnote{julialang.org} \citep{Bezanson_Julia_A_fresh_2017}
raison d'être is its combination of performance and
ease of use in the form of rapid prototyping.
This has been contrasted to efficient languages,
like Fortran, on one side, and accessible ones, like Python and Matlab, on the other.
Within the science community, there has been a shift towards accessible languages,
especially Python, for scripting.
Such languages are ideal for fast data analysis and light-weight calculations.
Writing such code in Fortran or C would require more time and cause some extra work
during the debugging phase, which would not be justified by any speed increase for small problems.

Julia has been lauded for delivering the advantages of both worlds \citep{Bezanson_Julia_A_fresh_2017}.
It is supposed to be as accessible as Python, while as fast as Fortran or C.
That makes it particularly attractive to scientists who
not only wish to perform simple calculations, but want to do
serious numerical simulations on large clusters using MPI
for the parallelization \citep{churavy2022bridging}.
This allows scientists to develop HPC code for large real-world simulations.

For process parallelization, MPI.jl \citep{byrne2021mpi} is a popular Julia library.
While it makes it possible to run on CPU clusters,
Julia code can also be run on different kinds of GPUs.
Using the CUDA.jl library \citep{besard2018juliagpu} we can make
use of Nvidia GPUs, while using the AMDGPU.jl
library\footnote{\url{https://github.com/JuliaGPU/AMDGPU.jl}}
we are not limited to one GPU vendor and can run on AMD GPUs too.
For efficient data storage, the
HDF5.jl\footnote{\url{https://github.com/JuliaIO/HDF5.jl}}
and ADIOS2.jl\footnote{\url{https://github.com/eschnett/ADIOS2.jl}}
libraries have been implemented in Julia as wrappers, giving it some of the necessary
tools to run on HPC clusters.

Benchmarks for Julia have largely focused on its serial performance
for scientific calculations (e.g.\ \citep{giordanokloeweretal2022}).
This is unfortunate, as more and more HPC clusters offer more comprehensive
support to run Julia code \citep{churavy2022bridging}.
\cite{Hunold_Steiner_2020} tested the MPI communication efficiency on up to \num{36} nodes
with 32 cores each, totaling \num{1152} CPU cores.
Scaling the message size, they showed that Julia's MPI implementation, in practice,
only causes very little overhead, as was also shown by \cite{giordanokloeweretal2022}.
\cite{Strauss-Bishnu-2023} solved the shallow water equations using a Julia and
a Fortran code.
For both codes they found similarly good scaling for up to \num{2048} cores
on a cluster.
Their Julia code has a factor two lower compute time, giving it a significant edge.
The closest performance and scaling tests to exascale on CPUs and GPUs
were performed by \cite{Godoy-Valero-Lara-2023}.
Unlike previous scaling tests, they used a real world problem,
a diffusion-reaction model in three dimensions.
They showed that their Julia CUDA/MPI code had little overhead
and scaled to \num{4096} GPUs on \num{512} nodes.

One numerical code that has been written in Julia is
Trixi.jl \citep{schlottkelakemper2020trixi, ranocha2022adaptive}.
It is a numerical simulation framework to solve conservation laws and non-conservative systems.
Currently, it can solve equations in 1, 2, and 3 dimensions
and on different kinds of grids.
Using MPI it can run on clusters on multiple nodes.
With its recent GPU implementations, it can also run on
clusters of GPUs.

To test Julia's abilities to run on exascale clusters, we
use the problem of the Taylor-Green vortex in
three dimensions, solving the compressible Euler equations.
As further test we simulate Alfv\'en waves in three dimensions,
solving the magnetohydrodynamic equations.
We run our test simulations on the three clusters
JURECA and JUWELS Booster at the Jülich Supercomputing Centre,
and MOGON NHR at the Johannes Gutenberg University Mainz.
For comparison, we run the same problems using the MPI paralellized
Fortran code FLUXO\footnote{\url{https://github.com/project-fluxo/fluxo}}
\citep{bohm2020entropy,rueda2021entropy}.

\section{Package loading}\label{sec:julia_compiler}

Julia uses a just-in-time (JIT) compiler built on top of the LLVM compiler infrastructure
\citep{lattner2004llvm}. Functions are compiled to machine code when they are initially
called, after which the compiled code is cached and reused for subsequent calls. Therefore,
Julia programs experience a certain startup latency associated with loading packages and
executing code for the first time.
Using custom system images, this latency can be significantly reduced.
System images store the loaded packages and compilation results for a selection of
packages and methods.
Julia ships with a system image that is used by default. It contains the compiler itself,
the Base library, several standard library packages, and compiled code to do common
operations like working with the REPL. However, it is also possible to create custom system
images with \texttt{PackageCompiler.jl}\footnote{\url{https://github.com/JuliaLang/PackageCompiler.jl}}.

A custom system image stores some or all packages of a given environment and, optionally,
all direct and indirect dependencies. Note that once the system image is created, the
versions of all packages in it are fixed and can only be changed by creating a new system
image. In addition, methods that should be compiled and stored in the system image can be
specified as well, by passing
a precompilation script to \texttt{PackageCompiler.jl}. This script is then executed, tracing all method
compilations, and later compiling and storing these exact methods in the custom system image.

In the HPC context, it is also possible to create a custom system image on a separate
machine before executing a job. If the CPU architectures of the involved machines are
different, the desired target CPU architecture can be passed to \texttt{PackageCompiler.jl}.

Custom system images are not only useful for reducing the compute time
due to compiling required methods beforehand: Since a system image consists of a single
large file, it can be easier to handle for the I/O system than loading many smaller
files that typical Julia packages are made of, depending on the exact parallel file system
and its configuration.

%%%%%%%%%%%%%%%%%%%%%%%%%%%%%%%%%%%%%%%%%%%%%%%%%%%%%%%%%%%%%%%%%%%%%%%%%%%%%%%%%%%%%%%%%%%%
%%%%%%%%%%%%%%%%%%%%%%%%%%%%%%%%%%%%%%%%%%%%%%%%%%%%%%%%%%%%%%%%%%%%%%%%%%%%%%%%%%%%%%%%%%%%
\section{Numerical simulation codes}\label{sec:simulation_codes}

For our scaling analysis we use two different numerical simulation packages, Trixi.jl and FLUXO.
Both solvers are designed for high-order discontinuous Galerkin spectral element
methods (DGSEM) applied to systems of conservative and non-conservative equations.
Additionally, Trixi.jl supports other DG schemes, such as those on triangular meshes.
Both tools employ low-storage Runge-Kutta (RK) methods for time integration.
They also support adaptive mesh refinement (AMR) using the \texttt{p4est}
library \citep{BursteddeWilcoxEtAl11} and use entropy- and energy--stable flux
formulations with subcell limiting strategies for shock capturing.
Both codes employ parallelization strategies based on MPI and, in the case of Trixi.jl,
an additional hybrid MPI-multithreading strategy is also available.
Below, we provide a more detailed description of each package.

In terms of single-thread performance, both codes are highly competitive.
Benchmark comparisons in curvilinear meshes indicate that they exhibit a similar
serial performance \citep{ranocha2021efficient} with a slight advantage
of Trixi.jl.

\subsection{Trixi.jl}
Trixi.jl \citep{schlottkelakemper2021purely, ranocha2022adaptive} is an open-source
Julia package for adaptive high-order simulations of conservation laws and
non-conservative systems.
It is used for computational fluid dynamics, including applications in
compressible fluid flow, astrophysics, aeroacoustics, and atmospheric flows.
A central focus of Trixi.jl is facilitating research on novel models and methods
while ensuring accessibility for new users and high performance for large--scale simulations.

Trixi.jl employs DGSEM for discretizing the governing equations on various mesh types,
including one- to three-dimensional meshes, Cartesian, curvilinear, and non-conforming grids,
as well as entropy-stable DG schemes for triangular meshes.
Time integration is performed using low-storage Runge-Kutta schemes.
More details on the numerical methods in Trixi.jl can be found in
\citep{schlottkelakemper2021purely, ranocha2021efficient}.

For parallel execution, Trixi.jl combines MPI and multithreading. 
The MPI implementation minimizes the number of communication operations,
using only one MPI communication operation at each right-hand-side computation
for hyperbolic conservation laws, at the cost of recomputing some numerical fluxes.
Non-blocking send and receive operations are used, and these communications are
overlapped with local computations to hide latencies.
The mesh is partitioned using a space-filling curve, ensuring a balanced load distribution.
A prototype for multi-GPU execution also exists, expanding Trixi.jl’s scalability
to heterogeneous computing environments.

\subsection{FLUXO}
FLUXO (see, e.g., \citep{bohm2020entropy,rueda2021entropy}) is a high-order
discontinuous Galerkin solver written in Fortran, designed for solving
both conservative and non-conservative equations.
Its primary applications include numerical astrophysics, space physics,
and compressible-flow aerodynamics, where accurate simulations of
complex physical phenomena are essential.
FLUXO is specifically designed for implementing entropy-stable and
kinetic-energy-preserving DG methods and is optimized for large-scale simulations.

% SC: Has FLUXO been tested on 'hundreds of thousands of cores'?
FLUXO operates solely on three-dimensional curvilinear and non-conforming meshes,
supporting AMR to enhance accuracy in critical regions and manage strong
shocks and discontinuities.
Its parallelization strategy, adapted from the massively parallel
Fortran code FLEXI \citep{krais2021flexi}, employs MPI without multithreading.
This approach enables efficient scaling to hundreds of thousands of cores while
maintaining a minimal computational load per rank \citep{beck2015high}.

Unlike Trixi.jl, FLUXO ensures that each numerical flux is computed only once,
with the effect that the computing load does not increase with the number of MPI ranks,
but requires two send and receive operations at each right-hand side computation
step for hyperbolic conservation laws.
This MPI implementation is colloquially called a "ping-pong" strategy,
as it requires neighboring ranks to communicate in one direction at the beginning
of the RHS computation and then again in the opposite direction after
the numerical fluxes are computed.
Non-blocking send and receive operations are used to reduce communication overhead.
To maintain balanced workloads, the domain is partitioned using a space-filling curve,
distributing shared interfaces evenly among ranks and minimizing MPI communication volume.

%%%%%%%%%%%%%%%%%%%%%%%%%%%%%%%%%%%%%%%%%%%%%%%%%%%%%%%%%%%%%%%%%%%%%%%%%%%%%%%%%%%%%%%%%%%%
%%%%%%%%%%%%%%%%%%%%%%%%%%%%%%%%%%%%%%%%%%%%%%%%%%%%%%%%%%%%%%%%%%%%%%%%%%%%%%%%%%%%%%%%%%%%
\section{Experimental setup}\label{sec:setup}

To evaluate the parallel performance of our code, we conduct scalability experiments using two three-dimensional test problems in periodic domains, incorporating both conservative and non-conservative equations.

\subsection{Taylor-Green vortex}

The first test problem involves the three-dimensional compressible Euler equations in conservative form:
\begin{eqnarray}
\partial_t
\begin{pmatrix} \rho \\ \rho \vec{v} \\ \rho e
\end{pmatrix}
+
\nabla \cdot
\begin{pmatrix} 
\rho \vec{v} \\
\rho (\vec{v}\, \vec{v}^{\,T}) + p\underline{I} 
 \\
\vec{v}\left(\frac{1}{2}\rho \left\|\vec{v}\right\|^2 + \frac{\gamma p}{\gamma -1}\right)
\end{pmatrix}
 = 
 \begin{pmatrix}
     0 \\ \vec{0} \\ 0
\end{pmatrix},
\end{eqnarray}
where $\rho$ is the fluid density, $\vec{v}$ the velocity, $e$ the specific total
energy and $p$ the pressure.
Pressure is not an independent variable and is related to the energy densities via
\begin{equation}\label{eq: pressure-energy}
p = (\gamma - 1) \left( \rho e - \frac{1}{2} \rho (v_x^2+v_y^2+v_z^2) \right),
\end{equation}
where $\gamma = 1.4$ is the adiabatic index.

As initial condition we choose the Taylor-Green vortex:
\begin{eqnarray}
\rho & = & 1 \\
v_x & = & \sin(x)\cos(y)\cos(z) \\
v_y & = & -\cos(x)\sin(y)\cos(z) \\
v_z & = & 0 \\
p & = & 100\rho/\gamma + 1/16 \rho(\cos(2x)\cos(2z) + 2\cos(2y) \nonumber \\
 & & + 2\cos(2x) + \cos(2y)\cos(2z)).
\end{eqnarray}
The conversion from $p$ to $e$ is done using equation \eqref{eq: pressure-energy}.
The computational domain is a periodic cube of size $(2 \pi)^3$ and origin at $(-\pi, -\pi, -\pi)$.

\subsection{Alfv\'en wave}
\label{sec: alfven wave}

The second test problem involves the three-dimensional compressible
magnetohydrodynamics (MHD) equations, which incorporate a generalized
Lagrange multiplier (GLM) divergence-cleaning technique
(see, e.g., \citep{derigs2018ideal, bohm2020entropy}).
These equations include non-conservative terms that arise when the
divergence of the magnetic field is nonzero at the discrete level.

The GLM-MHD equations extend the Euler equations by describing the evolution of the same
hydrodynamic quantities; density $\rho$, momentum components
($\rho v_x$, $\rho v_y$, $\rho v_z$), and specific total energy $e$;
while also accounting for the evolution of the magnetic field
$\vec{B} = (B_x, B_y, B_z)^T$ and the so-called divergence-cleaning field $\psi$.
The equations read as
\begin{strip}
\begin{eqnarray}\label{eq:glm-mhd}
\partial_t
\begin{pmatrix} \rho \\ \rho \vec{v} \\ \rho e
 \\ \blue{\vec{B}} 
 \\ \red{\psi} 
\end{pmatrix}
+
\nabla \cdot
%
%% Flux
%% --------------
\begin{pmatrix} 
\rho \vec{v} \\
\rho (\vec{v}\, \vec{v}^{\,T}) + p\underline{I} 
+ \blue{\frac{1}{2} \|\vec{B}\|^2 \underline{I} - \vec{B} \vec{B}^T}
 \\
\vec{v}\left(\frac{1}{2}\rho \left\|\vec{v}\right\|^2 + \frac{\gamma p}{\gamma -1}\right)
+ \blue{\left( \vec{v}\,\|\vec{B}\|^2 - \vec{B}\left(\vec{v}\cdot\vec{B}\right) \right)} 
+ \red{c_h \psi \vec{B}}  
 \\
 \blue{\vec{v}\,\vec{B}^T - \vec{B}\,\vec{v}^{\,T}} 
+ \red{c_h \psi \underline{I}} \\ 
 \red{c_h \vec{B}} \\ 
\end{pmatrix}
\nonumber\\
%% Powell source term
%% ------------------
 +  \blue{
(\nabla \cdot \vec{B}) 
\begin{pmatrix}
 0 \\ \vec{B} \\ \vec{v} \cdot \vec{B} \\  \vec{v} 
\\ \red{0}
\end{pmatrix}
 }
 +
 \red{
\nabla \psi \cdot 
\begin{pmatrix}
 0 \\ \underline{0} \\ \vec{v} \psi \\  \underline{0} \\ \vec{v}
\end{pmatrix}
}
 = 
 \begin{pmatrix}
     0 \\ \vec{0} \\ 0 \\ \vec{0} \\ 0
 \end{pmatrix},
\end{eqnarray}
\end{strip}
where $c_h$ is the divergence-cleaning speed, a parameter that is adjusted at
every time step to maximize the divergence-cleaning effect without restricting
the time-step size, $\underline{I}$ is the $3 \times 3$ identity matrix,
and the gas pressure is now computed as
\begin{eqnarray}\label{eq:pressure-mhd}
p & = & (\gamma - 1) \left( \rho e - \frac{1}{2} \rho (v_x^2+v_y^2+v_z^2) \right. \nonumber \\
 & & \left. \blue{
- \frac{1}{2} (B_x^2+B_y^2+B_z^2) }
\red{
- \frac{1}{2} \psi^2}
\right),
\end{eqnarray}
where $\gamma = 5/3$ is the adiabatic index for a monoatomic plasma.
In equations \eqref{eq:glm-mhd} and \eqref{eq:pressure-mhd},
we have colored in black the terms that come from the Euler equations,
in blue the terms coming from Maxwell's equations for electromagnetism,
and in red the terms that are only related to the GLM divergence-cleaning technique.

As initial condition, we choose the advection of an Alfvén wave:
\begin{eqnarray}
\rho & = & 1 \\
v_x & = & - 0.2 n_y \cos{\phi} \\
v_y & = & 0.2 n_x \cos{\phi}  \\
v_z & = & 0.2 \sin{\phi}  \\
p & = & 1 \\
B_x & = & n_x - v_x \\
B_y & = & n_y - v_y \\
B_z & = & - v_z \\
\psi & = & 0,
\end{eqnarray}
with the constants $n_x = 1 / \sqrt{5}$, $n_y = 2 / \sqrt{5}$, and
\begin{equation}
    \phi = \frac{2 \pi}{n_y} 
    \left( 
    n_x (x-1)
    +  n_y (y-1)
    \right).
\end{equation}
The computational domain is a periodic cube with size $2^3$ and origin at $(-1, -1, -1)$.

%%%%%%%%%%%%%%%%%%%%%%%%%%%%%%%%%%%%%%%%%%%%%%%%%%%%%%%%%%%%%%%%%%%%%%%%%%%%%%%%%%%%%%%%%%%%
%%%%%%%%%%%%%%%%%%%%%%%%%%%%%%%%%%%%%%%%%%%%%%%%%%%%%%%%%%%%%%%%%%%%%%%%%%%%%%%%%%%%%%%%%%%%
\section{Results}\label{sec:results}

\subsection{Taylor-Green vortex}

We evaluate the performance of the two codes through strong scaling tests where the problem size
is kept constant. The discretization is done via
a flux differencing DG scheme combined with an explicit 5-stage Runge-Kutta method for time
integration. The test problem runs for 200 time steps and uses \num{262144} elements
and a polynomial degree of 3 for the DG scheme, thus the total number of degrees of freedom is
$$
  \text{\#DOF} = \num{262144} \times (3 + 1)^3 = \num{16777216}.
$$
Parallelization is achieved via MPI without multithreading. For this type of
problem without mesh adaptivity, Trixi.jl does not allocate any additional memory during
the main loop and thus Julia's garbage collector is not triggered in the main loop.

The experiments were run on the JURECA and JUWELS Booster systems at the Jülich Supercomputing Centre (JSC),
and the MOGON NHR system at the Johannes Gutenberg University Mainz.
Their detailed architecture can be found in table \ref{table:hardware-configuration}.

\begin{table*}
  \centering
  \begin{tabularx}{\textwidth}{l X X X}
  \toprule
                    & JURECA                       & JUWELS Booster                   & MOGON NHR \\
  \midrule
       nodes        & 480                          & 936                              & 590 \\
       CPUs / node  & 2 × AMD EPYC 7742, 2.25 GHz  & 2 × AMD EPYC Rome 7402, 2.8 GHz  & 2 x AMD EPYC 7713, 2.0 GHz \\
       cores / node & 128                          & 48                               & 128 \\
       RAM / node   & 512 GB DDR4, 3200 MHz        & 512 GB DDR4, 3200 MHz            & 400 x 256GB, 159 x 512GB, 27 x 1024GB, 4 x 2048GB \\
       Interconnect & InfiniBand HDR100 (Connect-X6) & 4 x Infiniband HDR\linebreak (Connect-X6) & Infiniband HDR (Connect-X6) \\
  \bottomrule
  \end{tabularx}
  \caption{Hardware configuration of used systems.}
  \label{table:hardware-configuration}
\end{table*}

FLUXO was compiled with the Intel Fortran 2021.4 compiler with `-O3` optimizations enabled.
For the Trixi.jl tests we use two different library and Julia versions.
Runs on the JURECA cluster use OpenMPI v4.1.2 and Julia v1.8.3 with
bounds checking disabled, while runs on the JUWELS Booster and MOGON NHR
cluster use OpenMPI v4.1.5 and Julia v1.10.6 with bounds checking enabled.
All experiments use Trixi.jl v0.5.1.
For all the FLUXO simulations we used the commit \texttt{c7e0cc9b7fd4569dcab67bbb6e5a25c0a84859f1}.
Each simulation was run 5 times on both codes and the minimum runtime was used to
calculate performance metrics.
Since Julia compiles functions when they are called for the first
time, we discard the runtime of the first runs for Trixi.jl because we intend to compare
the pure runtimes without compilation times.

As main performance metric we use the number of degrees of freedom updates per second
(DOF updates/second or DUPS).
Other possible scaling metrics, like wallclock time or time per time step, would
not give us any additional information as they would be differently scaled or reciprocals.
The exception is the speedup, which compares the simulation runtime using a certain number
of MPI ranks with the runtime using the smallest amounts of ranks, i.e. ranks on a single node.
This can be 128 on the JURECA cluster, 48 on the JUWELS Booster machine and 128 on the MOGON NHR cluster.

\begin{figure}
  \centering
  \includegraphics[width=\linewidth]{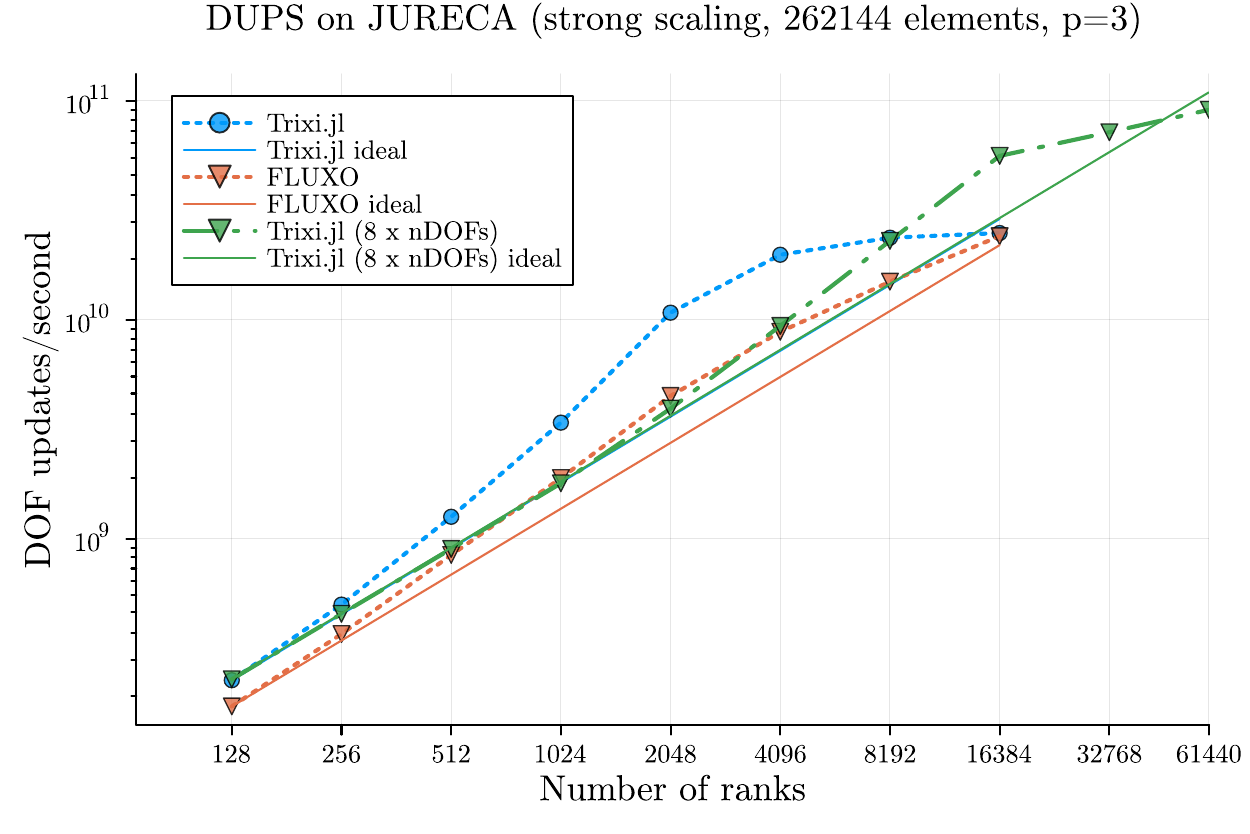}
  \caption{Degree of freedom updates per second for a flux differencing DG
    discretization of the 3D compressible Euler equations in FLUXO and Trixi.jl
    running on the JURECA cluster.
    Higher is better.}
  \label{fig:trixi-strong-scaling-jureca-large-dups}
\end{figure}

For the experiments on JURECA we show the scaling results for Trixi.jl and FLUXO
as DUPS in figure \ref{fig:trixi-strong-scaling-jureca-large-dups}.
Using the relatively small problem size we also get a good scaling behavior
up to \num{16384} ranks with an observed flattening off of the curve.
With a large number of MPI ranks communication overhead can be substantial,
which is why we also test an 8 times larger problem by increasing the initial
refinement level of the grid.
With that we obtain good scaling for the maximum testable number of MPI ranks
of \num{61440} for Trixi.jl, with a parallel efficiency of \num{0.83} relative to
the single-node baseline.
In this benchmark, Trixi.jl shows a superlinear speedup which is likely caused by cache effects.
Note that this superlinear behavior also affects the efficiency measured against the
single-node baseline: the efficiency peaks at \num{1.92} for \num{16384} ranks, and the
marginal efficiency beyond that point is only about \num{0.65} per doubling of the rank count.
Figure \ref{fig:trixi-strong-scaling-jureca-large-dups} shows that if the problem size is increased,
the superlinear bump is shifted accordingly to higher ranks.

\begin{figure}
  \centering
  \includegraphics[width=\linewidth]{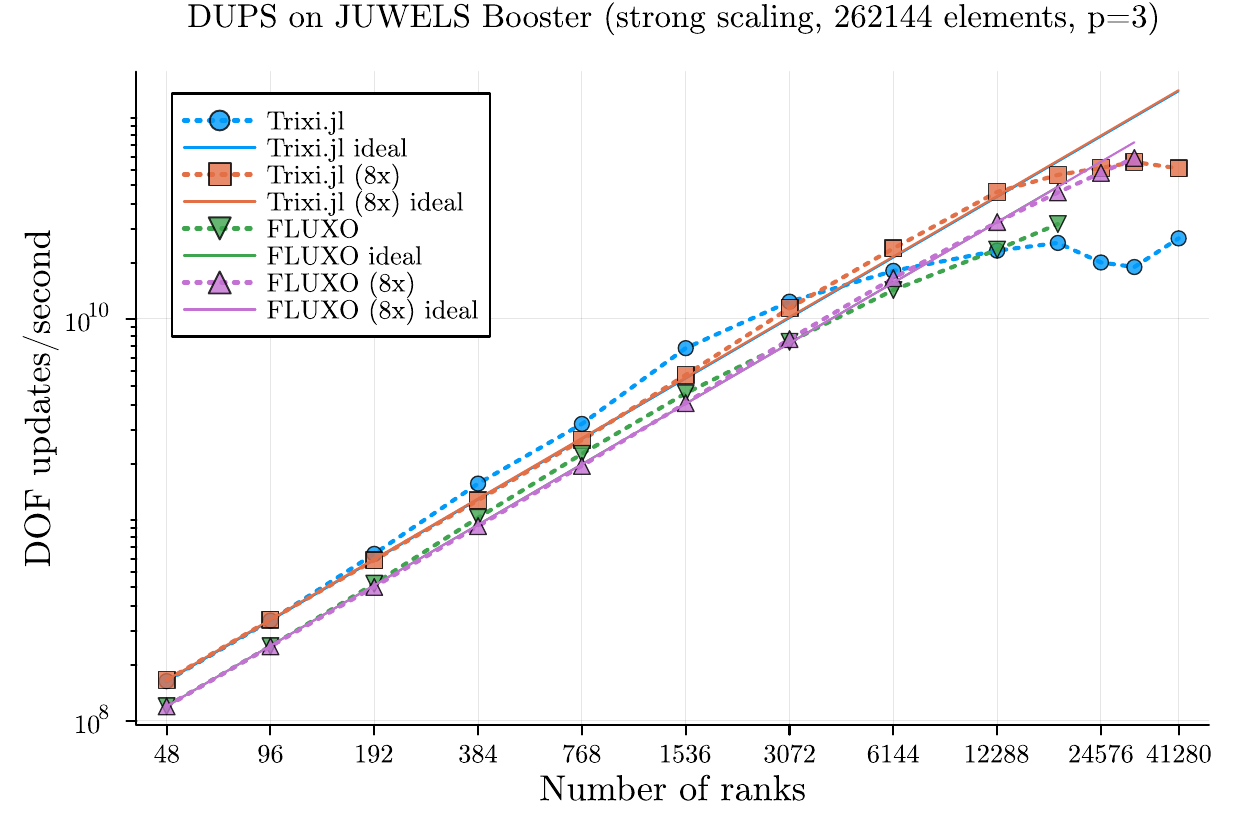}
  \caption{Degree of freedom updates per second for a flux differencing DG
    discretization of the 3D compressible Euler equations in FLUXO and Trixi.jl
    running on the JUWELS Booster cluster.
    Higher is better.}
  \label{fig:trixi-strong-scaling-juwelsbooster-dups}
\end{figure}

In figure \ref{fig:trixi-strong-scaling-juwelsbooster-dups} we show the DUPS for Trixi.jl and FLUXO on
the JUWELS Booster machine for the Taylor-Green vortex problem.
Both codes show good strong scalability from 48 MPI ranks (1 node) to \num{18432} MPI ranks (384 nodes).
For higher number of MPI ranks we see a clear flattening off, indicating poor scalability
beyond that number.
When the test problem size is increased by a factor of 8,
Trixi.jl shows good strong
scalability up to \num{18432} MPI ranks, retaining about \num{85}\% parallel
efficiency there, while beyond \num{30720} ranks (640 nodes) the throughput no longer
increases significantly, which we analyze in
section \ref{sec:dt-bottleneck}.
This shows the role that the problem size plays in the scalability.
Keeping the problem size low increases the overhead, like communication,
relative to the time spent on the calculation of the right hand side.

\begin{figure}
  \centering
  \includegraphics[width=\linewidth]{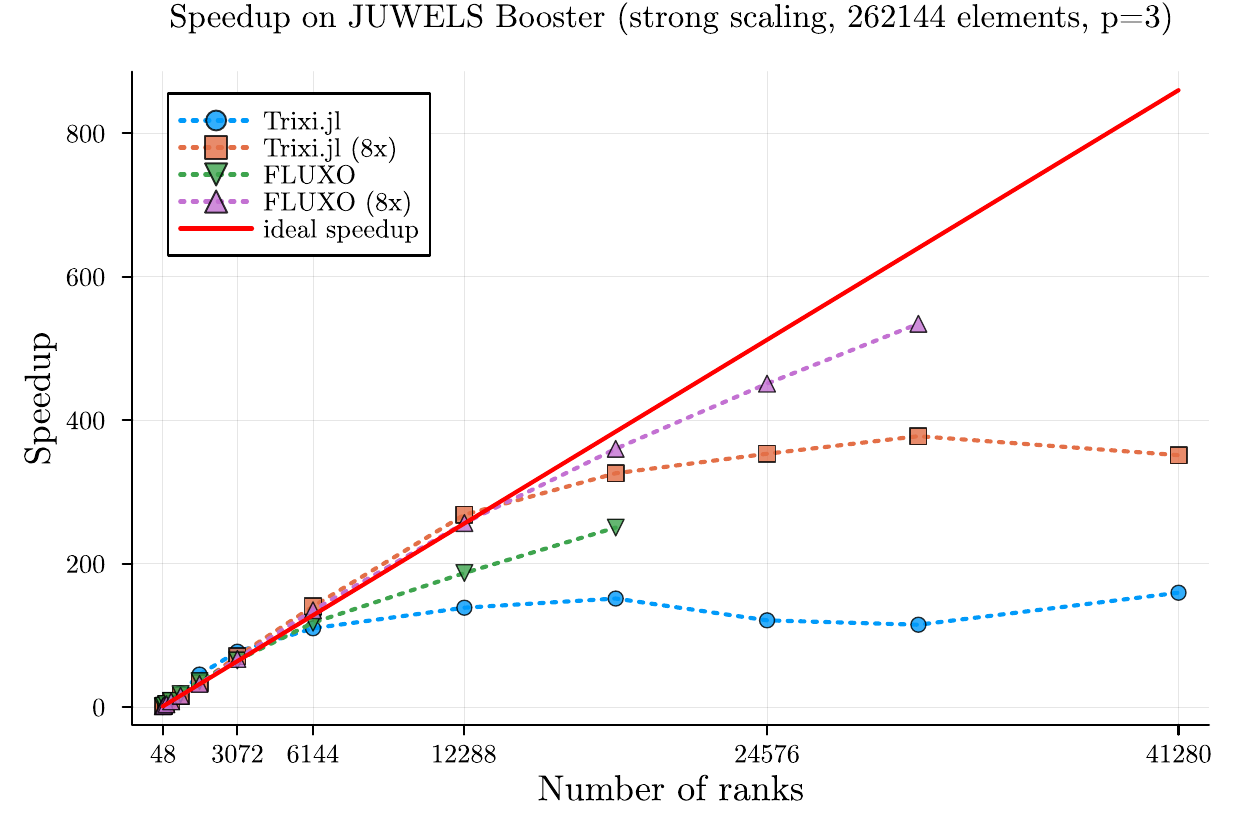}
  \caption{Speedup for a flux differencing DG discretization of the 3D compressible
    Euler equations in FLUXO and Trixi.jl running on the JUWELS Booster cluster.}
  \label{fig:juwelsbooster-spedup}
\end{figure}

Even though the overall runtime for Trixi.jl is smaller than for FLUXO, the scalability of FLUXO
seems to be slightly better.
Figure \ref{fig:juwelsbooster-spedup} shows that the speedup for
\num{18432} MPI ranks is higher for FLUXO.
However, this is unlikely to be related to Julia or Fortran,
but is most likely caused by the slightly faster serial performance of
Trixi.jl \citep{ranocha2021efficient} and by differences in the parallelization strategies of FLUXO and Trixi.jl.
The former employs a ping-pong strategy where two MPI processes communicate in one direction only,
then computations are done on the receiving rank only and the results get sent back.
Trixi.jl exchanges the data in both directions instead and does the calculations on each rank separately,
effectively performing them twice.
As the number of elements per rank gets smaller, this effect has
a larger impact on the total runtime.

\begin{figure}
  \centering
  \includegraphics[width=\linewidth]{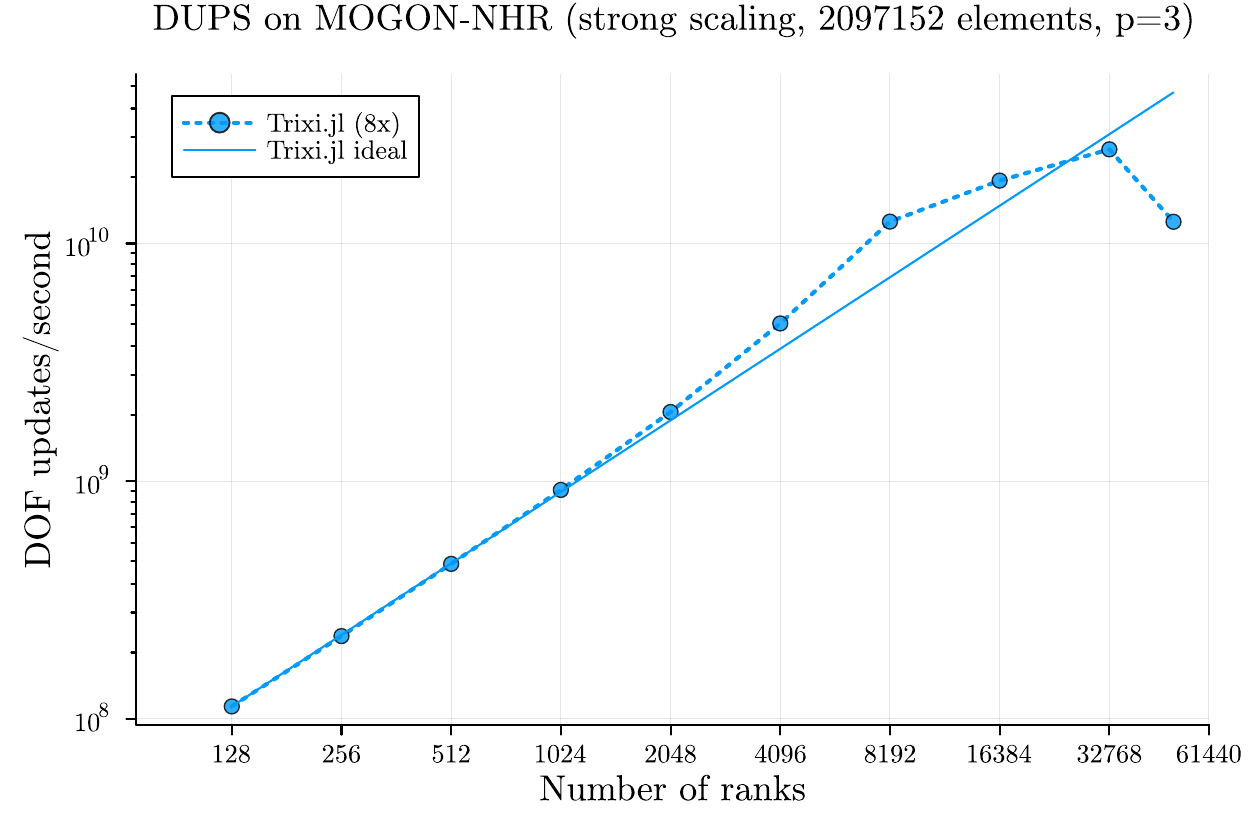}
  \caption{Degree of freedom updates per second for a flux differencing DG
    discretization of the 3D compressible Euler equations in Trixi.jl
     running on the MOGON-NHR cluster.
    Higher is better.}
  \label{fig:mogon-dups}
\end{figure}

\begin{figure}
  \centering
  \includegraphics[width=\linewidth]{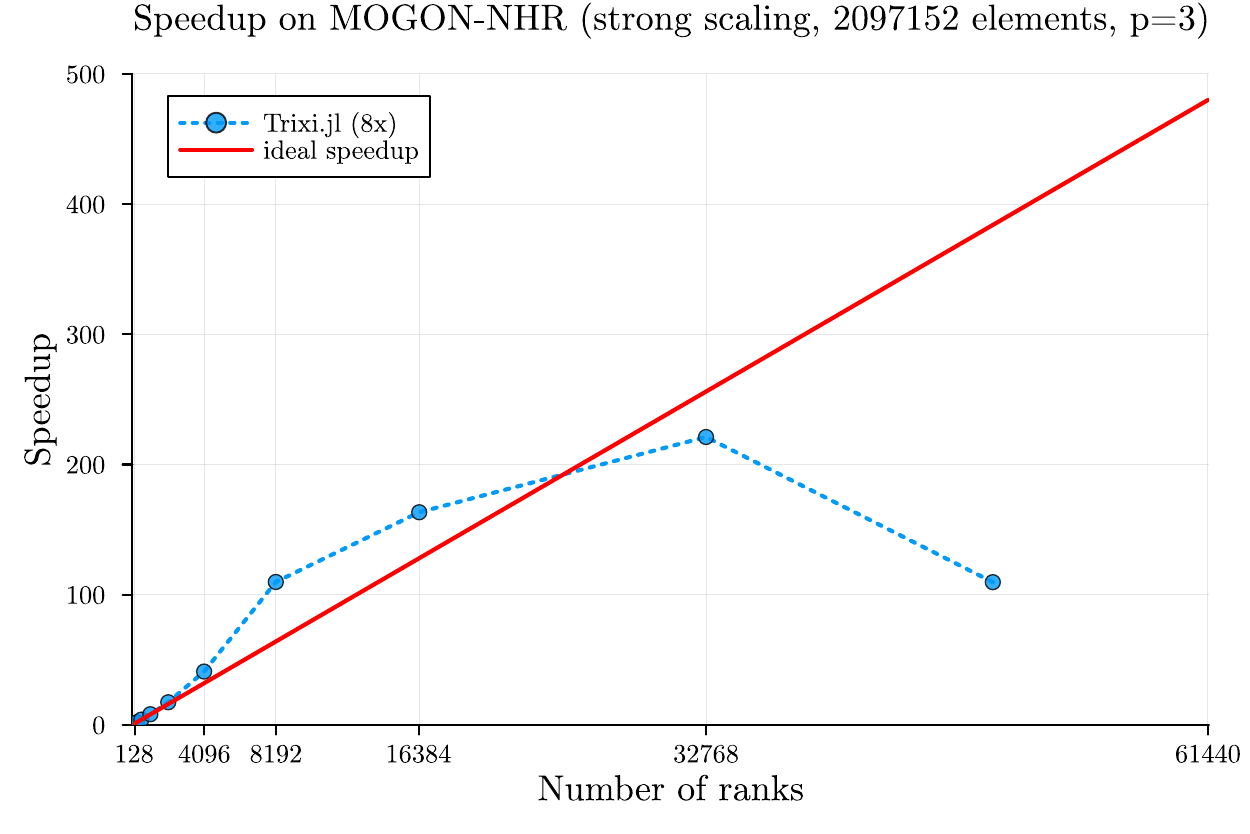}
  \caption{Speedup for a flux differencing DG discretization of the 3D compressible
    Euler equations in Trixi.jl running on the MOGON-NHR cluster.}
  \label{fig:mogon-spedup}
\end{figure}

Next, Figures \ref{fig:mogon-dups} and \ref{fig:mogon-spedup}
show the results of the experiments on the MOGON NHR cluster.
Here we run only the larger problem with $8\times$ the number of degrees of freedom.
Just as before, we get overall good scaling for up to \num{32768} MPI ranks;
beyond that the throughput drops.
This shows again that scaling can be significantly improved by increasing the problem
size per compute node.

\subsection{Alfv\'en wave}

Finally, we evaluate the HPC performance of Julia using the non-conservative problem
of the Alfv\'en wave in MHD.
Here we use the same resolution as for the larger Taylor-Green vortex problem
with \num{2 097 152} elements, which is eight fold larger
than the original problem.
With its large number of DOFs it already gives the problem a relatively large number
of calculations per rank.
Compared to the Euler equations from the Taylor-Green vortex test problem,
this MHD problem requires the solution of four additional variables.
Therefore, we solve a significantly more intense problem here, which should
show in an even better scaling.

\begin{figure}
  \centering
  \includegraphics[width=\linewidth]{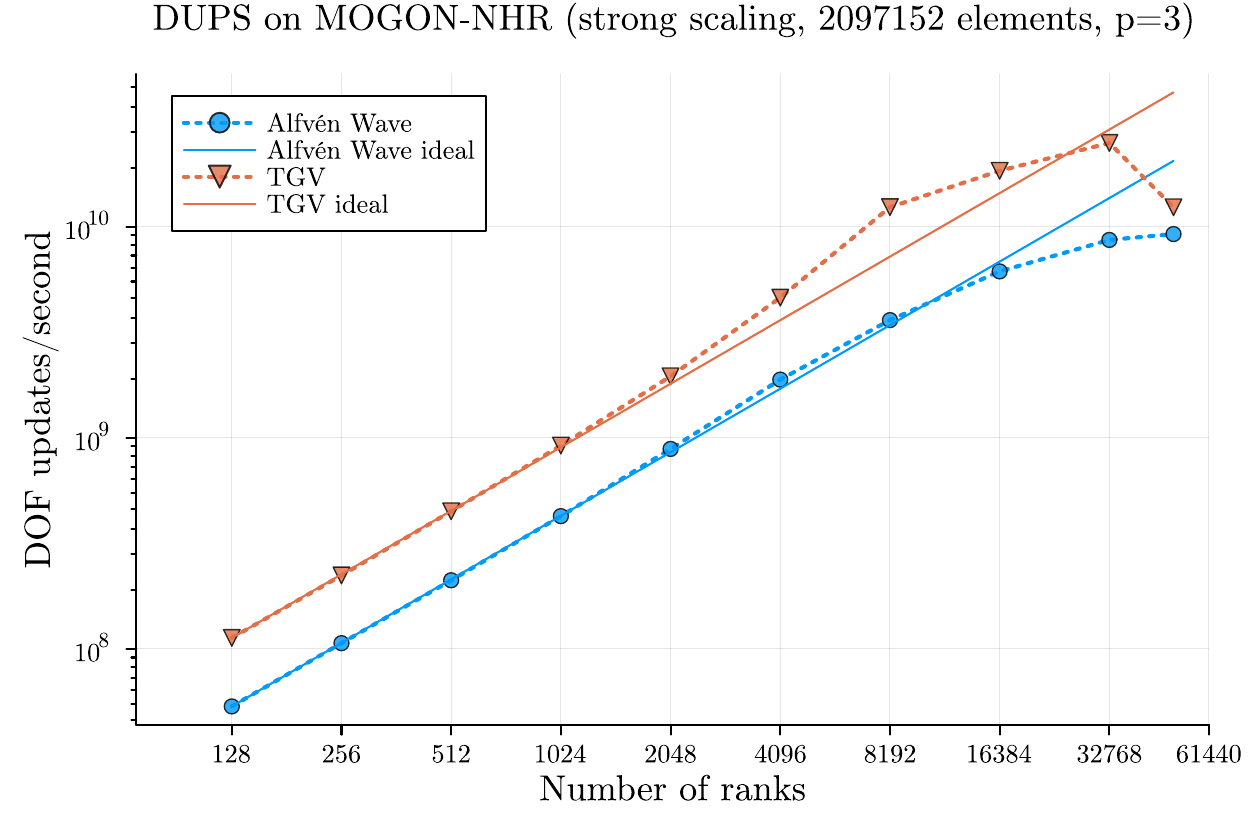}
  \caption{Degree of freedom updates per second for a flux differencing DG
    discretization of the 3D compressible MHD and Euler equations in Trixi.jl running on the MOGON-NHR cluster.
    Higher is better.}
  \label{fig:nc-mogon-dups}
\end{figure}

\begin{figure}
  \centering
  \includegraphics[width=\linewidth]{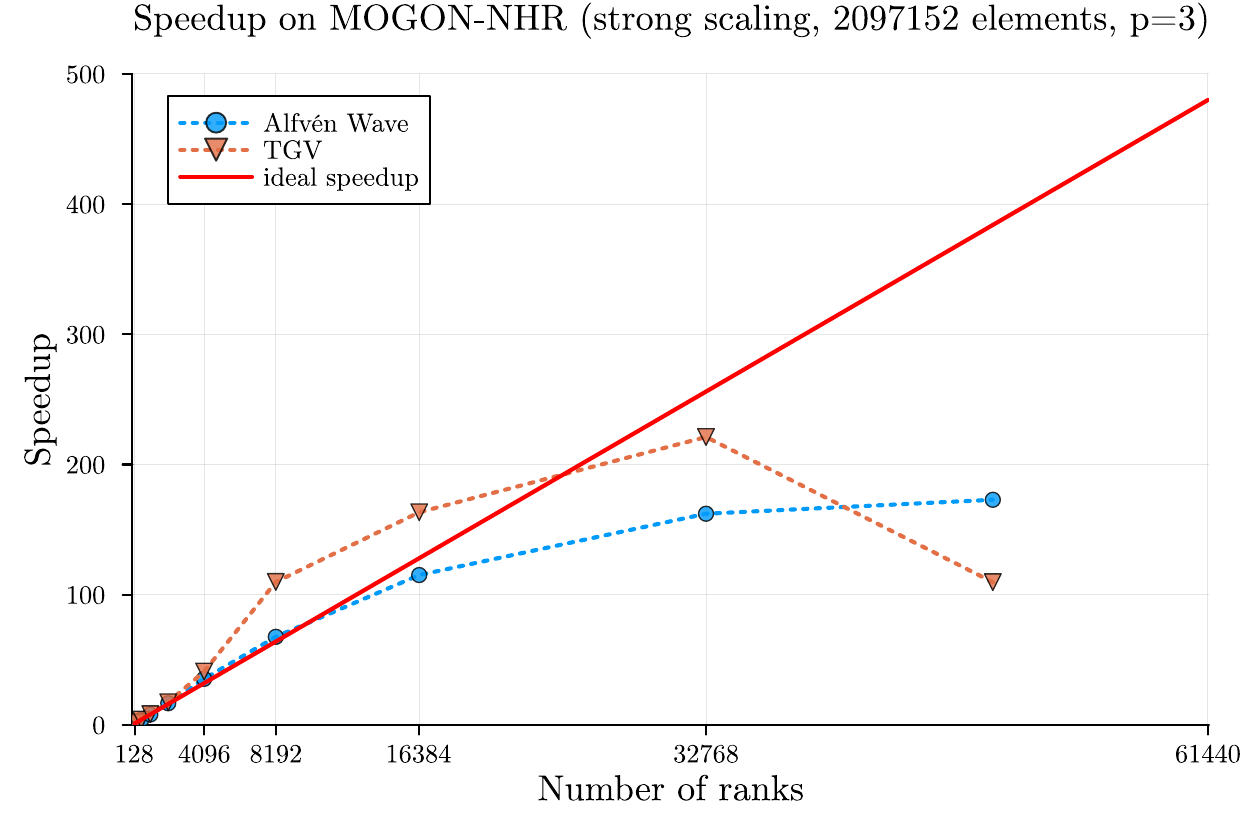}
  \caption{Speedup for a flux differencing DG discretization of the 3D compressible
    MHD equations in Trixi.jl running on the MOGON-NHR cluster.}
  \label{fig:nc-mogon-spedup}
\end{figure}

From Figures \ref{fig:nc-mogon-dups} and \ref{fig:nc-mogon-spedup} we see
an excellent strong scaling for this large problem up to \num{16384} ranks, with
throughput still increasing up to the maximum available
to us (\num{49152} ranks).

\subsection{Adaptive Time Step Bottleneck}\label{sec:dt-bottleneck}
Our scaling results show a clear bottleneck at high number of compute cores.
Depending on the problem solved and the cluster, good scaling is often limited to
a range of MPI ranks below the maximum we tested.
Larger problems scale better, as expected.
Communication overhead and rank boundaries become relatively smaller
compared to the computation of the right hand side.

So far the results shown have been for adaptive time stepping schemes.
Those require the regular calculation of the time step $\mathrm{d}t$ to fulfil the
CFL stability condition.
In Trixi.jl this happens for each rank and node through \texttt{StepsizeCallback}.
A call of \texttt{MPI.Allreduce} communicates the $\mathrm{d}t$ of each rank
and finds the global minimum, which is then used as the time step for all ranks.
This step should not require that much time, as MPI broadcasts
only a single float value.
However, we investigate the time required for this step and compare it
with the time required to calculate the right-hand side.

\begin{figure}
  \centering
  \includegraphics[width=\linewidth]{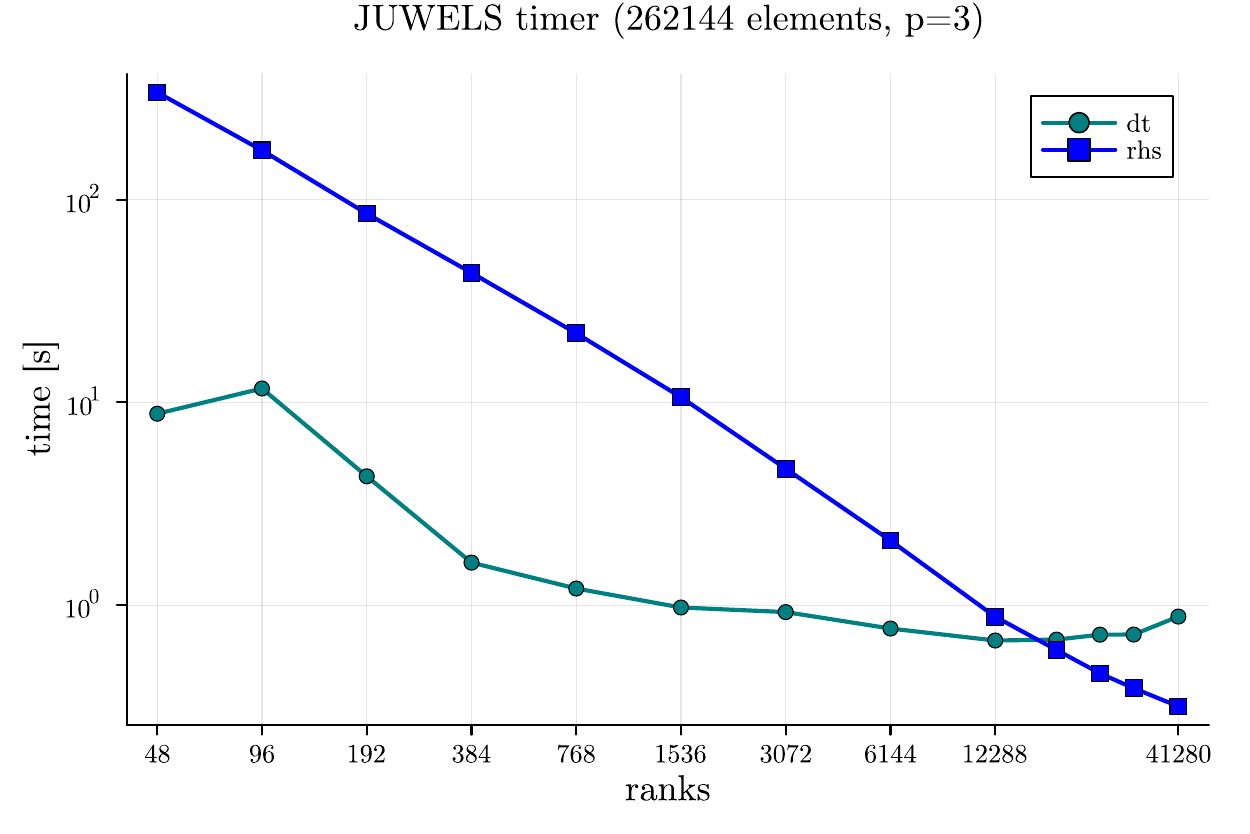}
  \caption{Timings for the computation of the time step and the right hand side on the JUWELS Booster cluster
  for the Euler test case.}
  \label{fig:trixi-strong-scaling-juwels-large-timer}
\end{figure}

Our timing tests are preformed using the Euler test case with \num{262144}
elements and polynomial degree 3.
On the JUWELS Booster cluster the time required for computing $\mathrm{d}t$ decreases with the number
of compute nodes until \num{2048} ranks (see Figure \ref{fig:trixi-strong-scaling-juwels-large-timer}).
Thereafter we see a relative flatting off with a small uptick at \num{41280} ranks.
At the same the time required to compute the RHS decreases steadily as expected.
The break-even occurs at ca.\ \num{18432} ranks.
If we compare this with Figure \ref{fig:trixi-strong-scaling-juwelsbooster-dups}
we see that at about the break-even point the scaling worsens.

\subsection{Fixed Time Step}
\begin{figure}
  \centering
  \includegraphics[width=\linewidth]{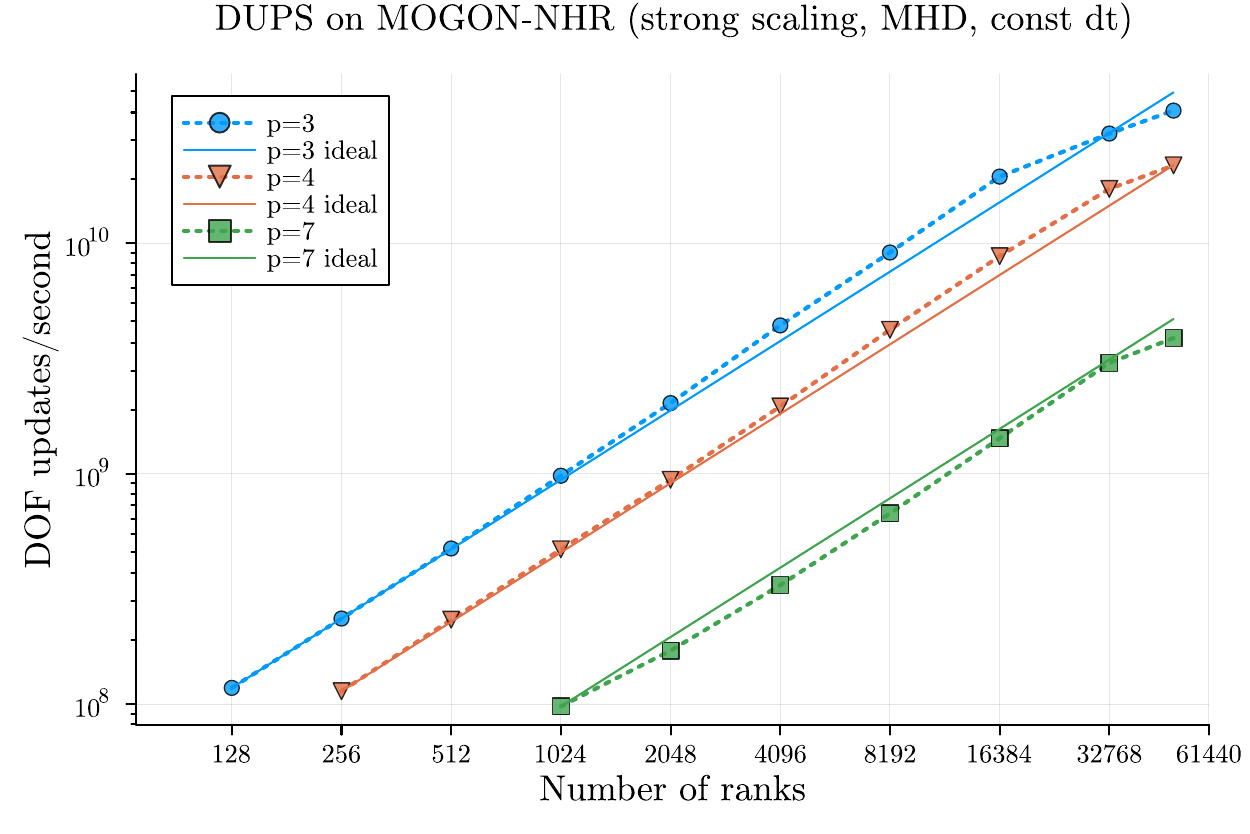}
  \caption{Degree of freedom updates per second for a flux differencing DG discretization of the 3D compressible MHD
  equations in Trixi.jl running on the MOGON-NHR cluster for fixed time step size $\mathrm{d}t$
  and different polynomial degrees $p$. Higher is better}
  \label{fig:trixi-strong-scaling-nc-mogon-large-mhd-dups}
\end{figure}

To test the scaling without the MPI bottleneck we run a series of benchmark calculations
with fixed time step $\mathrm{d}t$.
The test case is the Alfv\'en wave using polynomial degrees 3, 4 and 7 to increase
the problem size per time step.
Unlike the adaptive time step, we observe near ideal scaling for all core counts
available to us (see Figure \ref{fig:trixi-strong-scaling-nc-mogon-large-mhd-dups}).

\subsection{Avoiding the parallel I/O bottleneck during startup}
The code loading mechanism described in section \ref{sec:julia_compiler} poses challenges when
running a large number of Julia processes in parallel with MPI. The high-performance parallel
file systems usually deployed on supercomputers are generally not optimized for reading many
small files which is required during code loading and using a high number of processes may
lead to long loading times.
We investigate one possible way of avoiding the parallel I/O bottleneck during startup.

\subsubsection{Usage of a custom system image}
\begin{figure}
  \centering
  \includegraphics[width=\linewidth]{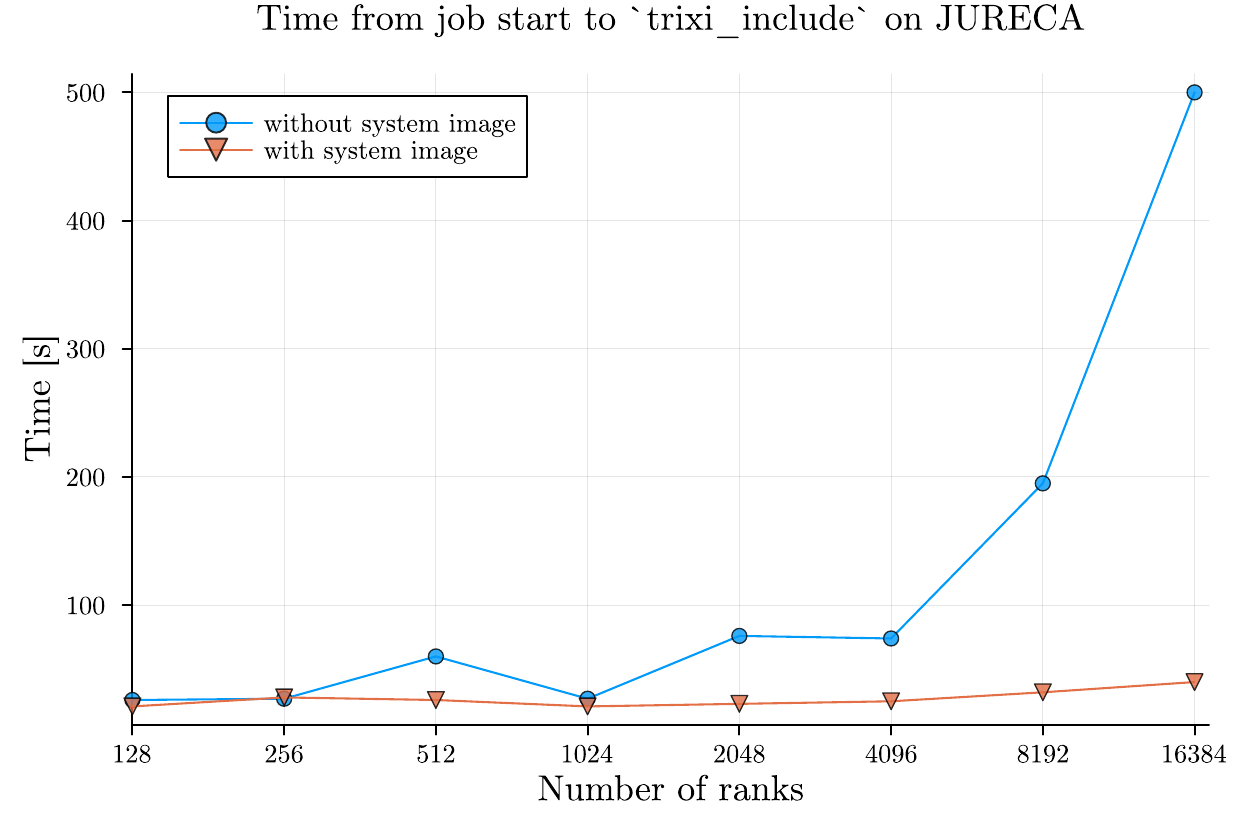}
  \caption{Time required for starting Julia and loading required packages for a
    flux differencing DG discretization of the 3D compressible Euler equations
    in Trixi.jl. For one setup a system image containing all packages and
    methods required for the simulation was used.}
  \label{fig:sysimg-time-until-include}
\end{figure}

The way of avoiding the parallel I/O bottleneck that we investigate is the use
of a custom system image.
As outlined in section \ref{sec:julia_compiler} it is possible
to create a custom system image that holds compiled versions of pre-specified methods.
In particular, it is possible to take a simulation setup for Trixi.jl and create a system
image that holds compiled versions of all methods that are required to run that particular
simulation.
The resulting system image is again a single large file which is easier to
handle for the parallel file system.
In addition to avoiding the I/O bottleneck, this method also reduces or, if done properly,
eliminates the compilation time, which is particularly
useful for smaller compute jobs where the compilation time makes up a larger portion of the
total runtime.
Figure \ref{fig:sysimg-time-until-include} shows that the package loading time
without a custom image starts to increase dramatically when the number of MPI processes is
increased but only slightly increases when a custom system image is used.

%%%%%%%%%%%%%%%%%%%%%%%%%%%%%%%%%%%%%%%%%%%%%%%%%%%%%%%%%%%%%%%%%%%%%%%%%%%%%%%%%%%%%%%%%%%%
%%%%%%%%%%%%%%%%%%%%%%%%%%%%%%%%%%%%%%%%%%%%%%%%%%%%%%%%%%%%%%%%%%%%%%%%%%%%%%%%%%%%%%%%%%%%
\section{Summary and conclusions}\label{sec:conclusions}
Running the Julia code Trixi.jl and Fortran code FLUXO on large clusters
we showed that Julia not only performs very well, but also scales well
for up to \num{61440} cores on 480 compute nodes using Julia's standard
MPI library.
Larger problems improve the scaling, as expected.
This was not free.
Parallel code loading in Julia with MPI causes a bottleneck that causes the
code loading time to rapidly increase at 8192 cores and beyond.
We fixed this issue by generating a system image that contains the pre-compiled
methods for the simulation.

%% file: acknowledgments.tex
\begin{itemize}

\item
Andrés Rueda-Ramírez gratefully acknowledges funding from the Spanish Ministry of Science, Innovation, and Universities through a ``Beatriz Galindo'' grant (BG23-00062), the German Federal Ministry of Education and Research (BMBF) under the WarmWorld Smarter program through the ``ICON-DG'' project (01LK2315B), and from the European Research Council through the Synergy Grant Agreement No.\ 101167322-TRANSDIFFUSE.

\item
Benedict Geihe's work has been funded by the German Federal Ministry of Education and Research (BMBF) within the project ADAPTEX under the funding code 16ME0668K,
and thereby also through the European Union -- NextGenerationEU. The views and opinions expressed are solely those of the author(s) and do not necessarily reflect
the views of the European Union or the European Commission. Neither the European Union nor the European Commission can be held responsible for them.

\item
Marco Artiano and Hendrik Ranocha were supported by the Deutsche Forschungsgemeinschaft (DFG, German Research Foundation, project numbers 513301895 and 528753982 as well as within the DFG priority program SPP 2410 with project number 526031774) and the Daimler und Benz Stiftung (Daimler and Benz foundation, project number 32-10/22). We acknowledge support from the Mainz Institute of Multiscale Modeling (M3ODEL) and from the Max Planck Graduate Center with the Johannes Gutenberg University of Mainz (MPGC).
        
\item
Simon Candelaresi has been supported by the Deutsche Forschungsgemeinschaft (DFG, German Research Foundation)
through the research unit FOR 5409 "Structure-Preserving Numerical Methods for Bulk- and Interface Coupling of
Heterogeneous Models (SNuBIC)" (project number 463312734).

\item
The authors gratefully acknowledge the Earth System Modelling Project (ESM) for funding this work by providing computing time on the ESM partition of the
supercomputer JUWELS \cite{JUWELS} at the Jülich Supercomputing Centre (JSC).

%Parts of this research were conducted using the supercomputer MOGON NHR and/or advisory services offered by Johannes Gutenberg University Mainz (hpc.uni-mainz.de), which is a member of the AHRP (Alliance for High Performance Computing in Rhineland Palatinate,  www.ahrp.info) and the Gauss Alliance e.V.
%The authors gratefully acknowledge the computing time granted on the supercomputer MOGON NHR at Johannes Gutenberg University Mainz (hpc.uni-mainz.de).

% \item
% The authors would like to thank the German Federal Ministry of Research, Technology and
% Space and the German federal states (http://www.nhr-verein.de/en/our-partners) for
% supporting this work/project as part of the National High-Performance Computing (NHR) joint
% funding program.
\item
The authors gratefully acknowledge the computing time made available to them on the high-performance computers MOGON NHR at the NHR Centers NHR Süd-West. These Centers are jointly supported by the Federal Ministry of Education and Research and the state governments participating in the NHR (\url{www.nhr-verein.de/unsere-partner}).
%
% \item Samuel Omlin (?)
% \item Valentin Sukera
% \item Gabriel Baraldi
% \item Joshua Lampert (?)
% \item Frederik Ekre (?)
\end{itemize}